\documentclass[12pt]{article}   	
\usepackage{geometry}                		
								\usepackage{amsmath}
\usepackage{amssymb}
\usepackage{amsthm}
\usepackage{amsfonts}

\def\la{\lambda}
\def\H{{\mathbb H}}
\theoremstyle{boldplain}

\newtheorem{cor}[equation]{Corollary}

\newtheorem{lemma}[equation]{Lemma}

\newtheorem{remark}[equation]{Remark}

\newtheorem{thm}[equation]{Theorem}
\newtheorem{theorem}[equation]{Theorem}
\newtheorem{defn}[equation]{Definition}

\def\al{\alpha}
\def\eps{\epsilon}
\def\R{\mathbb R}

\newenvironment{dedication}
  {      
   \thispagestyle{empty}
   
   \itshape             
   \raggedleft          
  }
  {\par 
  }

\begin{document}

\title{A note on laminations with symmetric leaves}
\author{Michael Kapovich}
\date{\today}	

\maketitle

\begin{dedication}
To Dennis Sullivan on the occasion of his 80th birthday, with great admiration.
\end{dedication}

\begin{abstract}
We prove that (apart from dimension $n=4$), each Riemannian solenoidal lamination with transitive homeomorphism group and leaves isometric to a symmetric space $X$ of noncompact type, is homeomorphic to the inverse limit of the system of finite covers of a compact locally-symmetric $n$-manifold. 
\end{abstract}

This note is motivated by a talk by Alberto Verjovsky on solenoidal manifolds in Cuernavaca in 2017 and the 
papers \cite{Sullivan-2014}, \cite{Verjovsky-2014}, \cite{Sullivan-Verjovsky} of Dennis Sullivan and Alberto Verjovsky. Our main result is Theorem \ref{thm:main} below describing $n$-dimensional homogeneous solenoidal laminations with leaves isometric to a symmetric space of noncompact type. A very detailed survey of solenoidal manifolds (especially in low dimensions) with a comprehensive bibliography of the subject, can be found in the paper \cite{Verjovsky} by Alberto Verjovsky. On a personal note, I also would like to add that many papers written by Dennis Sullivan on geometric topology, hyperbolic geometry, discrete groups and dynamics were an inspiration for much of my work ever since I was an undergraduate student in Novosibirsk. 

\medskip 
{\bf Acknowledgements.} This work was partially supported by the NSF grant DMS-16-04241. I am grateful to Steve Hurder for helpful references on geometry of solenoidal laminations and to the referees for their comments and suggestions.

\section{Generalities on laminations}

We refer the reader to \cite[section 1]{Candel} and \cite{Hurder}, especially, section 10, 
for in-depth treatment of laminations.  

An $n$-dimensional {\em topological} solenoidal lamination is a Hausdorff, 2-nd countable topological space $L$ equipped with a system of {\em local charts}, 
which are homeomorphisms $\phi_\al: U_\al\subset L \to V_\al\times C_\al\subset  \R^n \times C$, where  $C$ is the Cantor set, sets $C_\al$ ({\em local transversals}) are homeomorphic to $C$, $\{U_\al\}_{\al\in A}$ is an open cover of $L$, while each $V_\al\times C_\al$ is open in $\R^n \times C$. The preimages $\phi_\al^{-1}(V_\al\times \{c\})$ are 
the {\em local leaves} of $L$. 

We will be mostly interested in {\em compact} solenoidal laminations, but some of the discussion below is more general. 

For each lamination one defines an equivalence relation {\em generated} by the following (nontransitive) relation: $x\sim x'$ if and only if $x, x'$ belong to the same connected component of 
a common local leaf of $L$. 
Equivalence classes of this equivalence relation are the {\em leaves} of $L$; the leaf through $x$ is denoted $L_x$. 
Each leaf $L_x$ carries a natural topology with respect to which each leaf is a 
topological $n$-manifold. This topology is defined via its basis as follows. Pick $y\in L_x$ and a local chart 
$\phi_\al: U_\al\subset L \to V_\al\times C_\al$ such that $y\in U_\al$. There exists a unique $c\in C_\al$ such that 
$\phi_\al(y)\in V_\al\times \{c\}$. The product $V_\al\times \{c\}$ is  homeomorphic to $V_\al$ via the projection to the first factor. Then take a basis of the standard topology in $V_\al$, and take its preimage in $\phi_{\al}^{-1}(V_\al \times \{c\})$ under the composition
$$
 U_\al \stackrel{\phi_{\al}}{\longrightarrow} V_\al\times C_\al \to V_\al. 
$$
Doing so for all $y\in L_x$ and $U_\al$ containing such $y$ yields a basis of topology on $L_x$. 

A lamination is said to be {\em minimal} if each leaf is dense in $L$. 
A lamination is said to be {\em homogeneous} if its group of homeomorphisms acts transitively on $L$. 

\begin{remark}
More generally, one can consider topological laminations where $C$ is just a topological space and, accordingly, local transversals are not required to be totally disconnected, and instead of open subsets of $\R^n$ one takes open subsets of another model space. Then the transition maps are required to send subsets in $V_\al\times \{c\}$ to subsets of $V_\beta\times \{c'\}$. 
(In our setting, this property is automatic.) We will not discuss such general laminations in this note: By default, all laminations are assumed to be solenoidal. 
\end{remark}

In what follows, by a {\em map} of two laminations we will mean a continuous map $L\to L'$; since local transversals are totally disconnected, such a map necessarily  sends leaves to leaves. 

A {\em leaf-wise metric} on a lamination $L$ is a metrization of each leaf $L_x$ of $L$ by a metric $d_{L_x}$. For convenience, we extend these metrics to a metric, denoted $d$, on the entire $L$ by $d(x,x')=\infty$ if $x, x'$ belong to different leaves. The following definition appears to be nonstandard, I could not find it in the literature: 

\begin{defn}
A leaf-wise metric $d$ on $L$ is {\em continuous} if:  

\begin{enumerate}
\item Whenever two sequences $(x_i), (y_i)$ in $L$ satisfy $\lim_{i\to\infty} d(x_i, y_i)=0$, then for all sufficiently large $i$, the pairs of points $x_i, y_i$ belong to the same local transversal. 
(The local transversals, in general, will depend on $i$.) 

\item For any two sequences $(x_i), (y_i)$ which belong to the same $U_\al$ for all $i$ and such that $x_i, y_i$ belong to the same local transversals (depending on $i$) in $U_\al$, we have
$$
x_i\to x, y_i\to y \Rightarrow \lim_{i\to\infty} d(x_i,y_i)=d(x,y). 
$$
\end{enumerate} 
\end{defn}

In fact, all the leaf-wise metrics used in this paper will be not only continuous but also  {\em path-metrics} in the sense that the distance between any two points is the infimum of lengths of paths connecting these points. 

\medskip 
Similarly to topological laminations, one defines {\em smooth laminations}, by requiring that:

(a) The transition maps $\phi_\al\circ \phi_{\beta}^{-1}$ are smooth when restricted to each $\R^n\times \{c\}$.  

(b) All the partial derivatives of the transition maps (taken with respect to the $\R^n$-variables) are continuous as functions on $V_\beta\times C_\beta$. 

Accordingly, the leaves of a smooth lamination have natural structure of smooth $n$-manifolds. For smooth laminations one defines {\em Riemannian metrics} as leaf-wise Riemannian metrics which vary continuously with respect to the local transversals $C_\al$. (See \cite[section 1]{Candel}  for details.)  
If $L$ is compact, leaves $L_x$ have (uniformly, independently of $x$) bounded geometry. (Each point $x$ has a neighborhood in $L_x$ such that the metric on this neighborhood is uniformly bi-Lipschitz to the standard open unit ball in the Euclidean $n$-space. Moreover, the leafwise Riemannian curvature tensor has uniformly bounded norm.) In particular, the corresponding leaf-wise metric $d$ on $L$ is continuous.  By the partition of unity, every smooth solenoidal lamination admits a Riemannian metric. 

\medskip 
Below is an important class of solenoidal laminations. Start with a {\em compact} connected manifold $M$ and consider an infinite inverse system of nontrivial {\em finite regular} 
covering maps 
$$
... \to M_j\to  M_i\to M 
$$
such that each covering $M_i\to M$ corresponds to a finite index subgroup $G_i< \pi_1(M)$. (We use the linear order here is only for the notational convenience, the inverse system is actually arbitrary.) Let $\widehat{M}\to M$ denote the covering of $M$ corresponding to the intersection of  
all the subgroups $G_i$. Then the inverse limit of the above system of coverings is a 
compact connected solenoidal lamination $L$ with leaves homeomorphic to $\widehat{M}$. In the terminology of Hurder in \cite{Hurder}, $L$ is a {\em McCord lamination}: Such laminations were introduced and studied (in greater generality, the space $M$ was not assumed to be a manifold) by McCord in  \cite{McCord}. In particular, McCord proved that McCord laminations are minimal and homogeneous. Every McCord's lamination fibers over $M$ with totally disconnected fibers.

In the special case, when $G=\{1\}$, the manifold $\widehat{M}$ is the universal cover of $M$ and the lamination $L$ will be denoted $M_\infty$; I will refer to $M_\infty$ as a {\em McCord solenoid}. The existence of a system of finite-index subgroups in $\pi=\pi_1(M)$ whose intersection is trivial is known as the {\em residual finiteness property} of $\pi$.

The same construction works in the smooth setting, yielding a smooth solenoidal lamination. If $M$ had structure of a Riemannian manifold, we take pull-backs of the Riemannian metric to the covering spaces $M_i$. Accordingly, the inverse limit $L$ has structure of a Riemannian lamination: Each leaf of this lamination is a Riemannian covering space of $M$. In the case  of $L=M_\infty$, each leaf of $L$ is isometric to the universal Riemannian covering of $M$ and, hence, is quasiisometric to the fundamental group $\pi=\pi_1(M)$ (recall that $M$ is assumed to be compact). 

A (compact) topological manifold $M$, of course, need not be smoothable, hence, we cannot have a Riemannian metric. Nevertheless, as any Peano continuum, $M$ can be metrized via a {\em path-metric}\footnote{By compactness, this is equivalent to the property that the metric is geodesic.}: 
This result was conjectured (for general Peano continua) by Menger and proved independently by Bing \cite{Bing} and Moise \cite{Moise}. 
The length structure given by such a metric $d$ lifts to the covering spaces $M_i$, so that the maps $M_i\to M$ all preserve lengths of curves, i.e. are isometries of the length structures. 

Passing to the inverse limit, we obtain a leaf-wise path-metric $d_L$ on the corresponding McCord solenoid 
$L$, such that the restriction of the projection $L\to M$ to each leaf of $L$ is locally isometric, which implies that the metric $d_L$ is continuous. If $L=M_\infty$, the group $\pi$ acts isometrically on each leaf $L_x$ of $L$ and, since the metric on $L_x$ is a path-metric and the action is properly discontinuous and cocompact, $L_x$ is again quasiisometric to $\pi$. 

\medskip 
For completeness of the picture (even though, we will not need this), we note that one also has the notion of triangulated solenoidal laminations.  A {\em triangulation} of a solenoidal lamination $L$ is a triangulation of each leaf of $L$ so that simplices vary continuously with respect to the local transversals. 
 The following theorem is due to A.~Clark, S.~Hurder and O. Lukina, \cite{CHL}: 

\begin{thm}\label{thm:T1}
Every compact solenoidal Riemannian lamination $L$ admits a smooth triangulation of bounded geometry. 
\end{thm}

Here a triangulation of a Riemannian manifold is said to have bounded geometry if each $k$-simplex is $\la$-bilipschitz to the standard Euclidean $k$-simplex for some uniform constant $\la$.

Regarding the structure of solenoidal laminations, A.~Clark and S.~Hurder in  \cite[Theorem 1.2]{CH} proved:

\begin{thm}\label{thm:T2}
For every smooth compact  connected homogeneous solenoidal lamination $L$ with simply-connected leaves, there exists a compact topological manifold  $M$ such that the McCord solenoid $M_\infty$ is homeomorphic to $L$. 
\end{thm}

In addition to this difficult theorem we will need several easy properties of maps between laminations proven below.

\begin{lemma}\label{L1}
Suppose that $L, L'$ are two compact solenoidal laminations equipped with continuous leaf-wise metrics $d, d'$ respectively. Then every continuous map $f: L\to L'$ is uniformly continuous with respect to the metrics $d, d'$.
\end{lemma} 
\proof The proof is essentially a standard argument for uniform continuity of continuous functions on compact metric spaces. Suppose $f$ is not uniformly continuous. Then there exist two sequences $(x_i), (y_i)$ such that $d(x_i,y_i)\to 0$ (in particular,  for all sufficiently large $i$, 
$L_{x_i}=L_{y_i}$ and, moreover, by part 1 of the definition of a continuous metric, $x_i, y_i$ belong to the same local transversal), but $d'(f(x_i), f(y_i))>\eps$ for some  $\eps>0$ independent of $i$. By the compactness of $L$, we can assume that $x_i\to x, y_i\to y$ in the topology of $L$. Hence, by the continuity assumption on $d$ (part 2 of the definition), $x=y$. By the continuity assumption on $f$, $f(x_i)\to f(x), f(y_i)\to f(y)=f(x)$ and, moreover, for all large $i$, the pairs $x_i, y_i$ belong to common local transversals in the domain 
$U'_\al \subset L'$ of some local chart of $L'$. By applying the continuity assumption (part 2) for the metric $d'$, 
$$
\lim_{i\to\infty} d'(f(x_i), f(y_i))= d'(f(x), f(y))=0. 
$$
A contradiction. \qed 

\begin{cor}\label{C1}
Under the assumptions of Lemma \ref{L1}, every continuous map $f: L\to L'$ is {\em coarsely Lipschitz} in the sense that there exist constants $k\ge 1, a\ge 0$ such that for any two points $x, y\in L$, 
$$
d'(f(x), f(y))\le kd(x,y) + a. 
$$ 
\end{cor}
\proof Since $d, d'$ are path-metrics, this corollary is an immediate consequence of the uniform continuity of $f$, cf. \cite[Lemma 8.8]{Drutu-Kapovich}. \qed

\begin{cor}\label{C2}
Under the assumptions of Lemma \ref{L1}, every homeomorphism $f: L\to L'$ defines a quasiisometry between the leaves of $L, L'$ with respect to the metrics $d, d'$.
\end{cor}
\proof Since $f$ and $f^{-1}$ are coarse Lipschitz maps between the leaves of $L, L'$, the statement follows from one of the equivalent definitions of quasiisometries, see 
\cite[Definition 8.10]{Drutu-Kapovich}. \qed

\begin{lemma}\label{L2}
Suppose that $L$ is a compact Riemannian lamination and $L'=M_\infty$ is a topological McCord solenoid, homeomorphic to $L$. Then each leaf of $L$ (equipped with the leaf-wise Riemannian distance function $d$) is quasiisometric to the fundamental group of $M$ (equipped with the word-metric). 
\end{lemma}
\proof We equip $L'=M_\infty$ with the leaf-wise path-metric $d'=d_{M_\infty}$  defined earlier. Each leaf $L'_{x'}$ of $L'$   with the metric $d'$ is quasiisometric to $\pi_1(M)$. By Corollary 
\ref{C2}, every leaf of $L$ is quasiisometric to its image leaf in $L'$ (equipped with the metric $d'$). Lemma follows. \qed 

\section{Homogeneous solenoidal laminations with symmetric leaves} 

In this section we prove the main result of the note. We fix $X$, a symmetric space of noncompact type.  

\begin{theorem}\label{thm:main}
Let $L$ be a compact homogeneous solenoidal $n$-dimensional Riemannian lamination with  leaves isometric to $X$. Then there exists a closed aspherical  
$n$-manifold $M$ such that:

1. The McCord solenoid   $M_\infty$ is homeomorphic to $L$. 

2. If $n\ne 4$ then $M$ is homeomorphic to the quotient of $X$ by a discrete, torsion-free cocompact group of isometries. If $n=4$, then $M$ is homotopy-equivalent to such a quotient. For example, if $X=\H^n$ is the hyperbolic $n$-space, then $M$ is either homeomorphic (in dimensions $n\ne 4$) or homotopy-equivalent (in dimension $n=4$) to a hyperbolic $n$-manifold.  
\end{theorem}
\proof As the reader will observe, my contribution to this result is minimal, it is mostly to combine deep work by others and to prove Lemma \ref{L2}.  

Part 1 of the theorem is due to Clark and Hurder, see Theorem 1 above. 

We now proceed to Part 2. The manifold $M$ satisfies two properties: It is aspherical (since its universal covering space is homeomorphic to $X$ which is aspherical), in particular, 
its fundamental group $\pi=\pi_1(M)$ is torsion-free, and the group $\pi$ is quasiisometric to the leaves of $L$, i.e. to $X$ (see Lemma \ref{L2}). 

It now follows from quasiisometric rigidity of symmetric spaces of noncompact type that $\pi$ is isomorphic to a uniform lattice in the isometry group of $X$. This deep result is a combination of work of many people: 

\begin{itemize}
\item Tukia, \cite{Tukia88},  Gabai, \cite{Gabai}, Casson and Jungreis, \cite{CJ} for 
$X=\H^2$, Sullivan, \cite{Sullivan}, for $\H^3$, Tukia, \cite{Tukia86}, for $X=\H^n$ (see also \cite{Kapovich14} and \cite[Chapter 23]{Drutu-Kapovich}).  

\item R.~Chow, \cite{Chow}, for complex-hyperbolic spaces $X={\mathbb C} \H^n$.

\item P.~Pansu, \cite{Pansu}, for quaternionic-hyperbolic spaces and octonionic hyperbolic plane. 

\item B.~Kleiner and B. Leeb, \cite{KL1}, when $X$ is a symmetric space of higher rank without rank 1 factors. A bit later, this result was also obtained by A. Eskin and B. Farb in \cite{EF} 
by a different method.

\item Lastly, B.~Kleiner and B. Leeb, \cite{KL2}, handled the case of symmetric spaces with rank one factors. 
\end{itemize} 

We also refer the reader to surveys of these results in \cite{Farb} and  \cite[Chapter 25]{Drutu-Kapovich}. 

\medskip
Thus, we conclude that $\pi$  is isomorphic to the fundamental group of a closed locally-symmetric (modeled on $X$) 
$n$-manifold $M'$. In particular, $M$ is homotopy-equivalent to $M'$. In dimensions different from $4$, a compact $n$-dimensional manifold homotopy-equivalent to a locally-symmetric manifold of nonpositive curvature is actually homeomorphic to it: In dimension 2 it is classical, in dimension 3 it is a corollary of Perelman's Geometrization Theorem, and in dimensions $>4$ this result is due to Farrell and Jones, \cite{FJ1, FJ2}. \qed

\noindent Address: Department of Mathematics, \\
University of California, Davis\\
CA 95616, USA\\
email: kapovich@math.ucdavis.edu

\end{document}